# ROBUST NEAREST-NEIGHBOR METHODS FOR CLASSIFYING HIGH-DIMENSIONAL DATA

By Yao-ban Chan and Peter Hall

*University of Melbourne*

We suggest a robust nearest-neighbor approach to classifying high-dimensional data. The method enhances sensitivity by employing a threshold and truncates to a sequence of zeros and ones in order to reduce the deleterious impact of heavy-tailed data. Empirical rules are suggested for choosing the threshold. They require the bare minimum of data; only one data vector is needed from each population. Theoretical and numerical aspects of performance are explored, paying particular attention to the impacts of correlation and heterogeneity among data components. On the theoretical side, it is shown that our truncated, thresholded, nearest-neighbor classifier enjoys the same classification boundary as more conventional, non-robust approaches, which require finite moments in order to achieve good performance. In particular, the greater robustness of our approach does not come at the price of reduced effectiveness. Moreover, when both training sample sizes equal 1, our new method can have performance equal to that of optimal classifiers that require independent and identically distributed data with known marginal distributions; yet, our classifier does not itself need conditions of this type.

**1. Introduction.** In classification problems where sample size is much smaller than dimension, nearest-neighbor methods, after truncation to reduce noise, can enjoy particularly good performance. They have the potential to be highly adaptive, not least because they do not require explicit assumptions about marginal distributions.

However, in very high-dimensional settings, conventional nearest-neighbor methods can be adversely affected by "noise" from vector components that do not carry useful information for classification. Moreover, they are not robust against outliers. In particular, they can be influenced considerably by heavy-tailed features of sampling distributions and can fail to give accurate classification when marginal distributions do not enjoy finite variance. Their sensitivity to correlation among data components, particularly in very high-dimensional contexts, is not well understood. And their performance in high-dimensional, highly heterogeneous cases, where the tails of distributions can vary from very light to very heavy within the same data vector, is largely unknown.

These phenomena occur often in the area of gene microarray analysis. Each microarray represents thousands of gene expression levels, but the sample size is typically small. Furthermore, the underlying distributions of the gene expressions levels are generally unknown and are likely to be heterogeneous, heavy-tailed and significantly dependent upon each other. With these features, conventional nearest-neighbor methods for analysis are likely to be ineffective.

In this paper, we shall suggest a robust nearest-neighbor classifier, where thresholding and truncation to zeros and ones are used to increase performance and, in particular, to remove sensitivity to heavy-tailed behavior. Choosing the threshold appropriately is the key to good classification accuracy. Threshold selection must adapt both to distribution type and to the ways in which populations differ from one another. We suggest a simple and practicable approach to selecting the threshold. Unlike cross-validation, our technique gives good performance even when there is only one training data-vector from each population.

We shall use theoretical arguments and numerical simulation to show that our technique is relatively insensitive to dependence among vector components, and that it enjoys good classification accuracy in high-dimensional, highly heterogeneous cases. In settings such as these, the performance of truncated nearest-neighbor classifiers can surpass that of competitors, such as methods based on extrema or on false-discovery rate (FDR) ideas. The latter two approaches are often identical; see Jin [14] and Donoho and Jin [5].

Nearest-neighbor methods are popular because of the wide variety of data types for which they are appropriate. Their implementation requires only a measure of distance and, in particular, is not founded on distributional properties of the data. Therefore, nearest-neighbor classifiers enjoy a high degree of acceptance in settings involving complex data, for example, in pattern recognition. See Dasarathy [3] and Shakhnarovich, Darrell and Indyk [17], for instance.

Properties of nearest-neighbor classifiers in classical settings, where dimension is small relative to sample size, are quite well understood. See, in particular, Devroye, Györfi and Lugosi [4]. Chapter 5 of that monograph is an excellent guide to the literature. There is a very large number of papers on nearest-neighbor methods in other settings, and it is possible to mention only a few of them here. Early contributions include those of Cover and Hart [2] and Cover [1], who gave upper bounds to risk, and Wagner [18] and Fritz [7], who derived convergence properties of the error rate. Psaltis, Snapp and Venkatesh [16] extended Cover's [1] results to higher dimensional settings, but still with the number of dimensions much less than sample size. Kulkarni and Posner [15] and Holst and Irle [10] discussed the case of dependent data vectors.

In these relatively classical treatments, it is common to regard the order, $k$, of a nearest-neighbor classifier as a tuning parameter and, perhaps, to attempt to optimize over it. However, in a variety of contemporary applications the number of data in each sample is so small, especially relative to dimension, that there is little point in taking $k$ larger than 1. We argue that, in such cases, the information that is critical to good performance is accumulated not through the number of data, but through the many components of each data vector. With that in mind, in this paper we shall optimize performance in a way that is sensitive to dimension, rather than to sample size.

## 2. Methodology.

2.1. *Sparsity and truncation.* Assume we observe random $p$-vectors $X_1, \ldots, X_m$ and $Y_1, \ldots, Y_n$, drawn from $X$- and $Y$-populations, respectively. We wish to construct a classifier, for the purpose of ascribing a new $p$-vector, $Z$ say, to either population.

Suppose it is known that the respective components of $X$ and $Y$ distributions are similar, except that one of them has, for a potentially sparsely arrayed sequence of component indices, generally higher mean than the other. We can formalize at least part of this assumption, by asking that, if $X = (X^{(1)}, \ldots, X^{(p)})^\mathrm{T}$ and $Y = (Y^{(1)}, \ldots, Y^{(p)})^\mathrm{T}$, then,

(2.1) for each $k$, (a) $X^{(k)} - E(X^{(k)})$ and $Y^{(k)} - E(Y^{(k)})$ have similar distributions, and (b) $E(Y^{(k)}) \geq E(X^{(k)})$; and, for a potentially sparsely distributed sequence of indices $k$, (c) $E(Y^{(k)}) > E(X^{(k)})$.

The one-sided nature of parts (b) and (c) of (2.1) motivates a one-sided classifier. Alternative methodology and theory, very similar to that which we shall develop below, are available in the two-sided case.

In view of the possible sparsity, it seems reasonable to truncate components of the data vectors by deleting those that do not attain a threshold, $t$ say. This has the effect of reducing the amount of noise that is present in coordinate values that convey little or no information for classification.

There are a variety of ways of implementing a procedure such as this. For example, we may do it by replacing each component by 0 if it is less than $t$, and by 1 otherwise. That is, defining $I_i^{(k)} = I(X_i^{(k)} > t)$, $J_i^{(k)} = I(Y_i^{(k)} > t)$, $I_i = (I_i^{(1)}, \ldots, I_i^{(p)})^{\mathrm{T}}$ and $J_i = (J_i^{(1)}, \ldots, J_i^{(p)})^{\mathrm{T}}$, we may build the classifier using the indicator data vectors $I_1, \ldots, I_m, J_1, \ldots, J_n$. Alternatively, we could base the classifier on $X_1', \ldots, X_m', Y_1', \ldots, Y_n'$, where $X_i^{(k)\prime} = X_i^{(k)} I(X_i^{(k)} > t)$, $Y_i^{(k)\prime} = Y_i^{(k)} I(Y_i^{(k)} > t)$, $X_i' = (X_i^{(1)\prime}, \ldots, X_i^{(p)\prime})^{\mathrm{T}}$ and $Y_i' = (Y_i^{(1)\prime}, \ldots, Y_i^{(p)\prime})^{\mathrm{T}}$.

Using the indicator data, we could conclude that $Z$ came from the $X$ population if

$$(2.2) \quad \min_{1 \leq i \leq m} \sum_{k=1}^{p} (I_i^{(k)} - K^{(k)})^2 \leq \min_{1 \leq i \leq n} \sum_{k=1}^{p} (J_i^{(k)} - K^{(k)})^2,$$

where $K^{(k)} = I(Z^{(k)} > t)$; and that $Z$ was from the $Y$ population otherwise. Alternatively, in place of (2.2) we could use the criterion

$$(2.3) \quad \min_{1 \leq i \leq n} \sum_{k=1}^{p} (X_i^{(k)\prime} - Z^{(k)\prime})^2 \leq \min_{1 \leq i \leq m} \sum_{k=1}^{p} (Y_i^{(k)\prime} - Z^{(k)\prime})^2,$$

where $Z^{(k)\prime} = Z^{(k)} I(Z^{(k)} > t)$. In this case, if (2.3) were true, then we would conclude that $Z$ was from the $X$ population. However, relative to methods based on (2.2), this approach would suffer more from stochastic variability and, hence, be less robust, in cases where $X$ and $Y$ had heavy-tailed distributions.

2.2. *Empirical choice of $t$.* We suggest a method based on thresholding, as follows. Let $i_X$ and $i_Y$ denote the respective values of $i$ at which the minima on the left- and right-hand sides of (2.2) are achieved. In this notation, (2.2) is equivalent to $T \leq 0$, where

$$(2.4) \quad T = T(t) = \sum_{k=1}^{p} (I_{i_X}^{(k)} - J_{i_Y}^{(k)})(1 - 2K^{(k)}).$$

Let $\xi_p$ denote a sequence diverging to infinity; put

$$(2.5) \quad z_p = \xi_p \log p,$$

denoting a threshold; let

$$(2.6) \quad S^2 = S(t)^2 = \sum_{k=1}^{p} (I_{i_X}^{(k)} + J_{i_Y}^{(k)})$$

and define $t = \theta$ by:

(2.7)   $\theta$ is the infimum of values $t \geq 0$ such that $|T(t)|/S(t) > z_p$; or, if no such $t$ exists, take $t$ to be a default value, for example, $t = 0$ or $t = -\infty$.

In the case of independent components, it is feasible to use a smaller threshold, defining $z_p$ by

$$z_p = \xi_p (\log p)^{1/2}, \tag{2.8}$$

where $\xi_p \to \infty$, instead of by (2.5). Nevertheless, (2.5) is also appropriate in the case of independence. If, in (2.7), it were necessary to pass to the default value, then we would conclude that the classification problem was marginal. That is, there was insufficient information to solve the problem reliably.

With $t = \theta$ given by (2.7), the classifier suggested by (2.2) is as follows:

(2.9) Classify $Z$ as coming from the $X$ population, if $T(\theta) \leq 0$ and as coming from the $Y$ population otherwise.

Our theoretical justification for (2.5) will be based on the assumption that the components $X^{(k)}$ and $Y^{(k)}$ are produced by a generalized form of an infinite-order moving average. The generalization permits marginal distributions to vary extensively from one component to the next, so that they are heavy-tailed for some indices $k$ but light-tailed for others. Alternative models for weak dependence, for example based on autoregressive processes, can be shown to also lead to the threshold choice at (2.5).

From at least a theoretical viewpoint, exact choice of $\xi_p$ is largely unimportant. Any sequence, for example, $\xi_p = \log p$, which diverges more slowly than any polynomial is appropriate. In this way, the sensitivity of the tuning-parameter selection problem is greatly reduced; we pass from the parameter $t$, to which the classifier is very sensitive, to $\xi_p$, to which the classifier is largely insensitive. Practical, empirical choices of $\xi_p$ will be discussed in Section 4.

Motivation for a threshold-based approach to choosing $t$ can be provided as follows. Neglecting, for the moment, the fact that $i_X$ and $i_Y$ at (2.4) are random variables; taking the components $I_{i_X}^{(k)}$ and $J_{i_Y}^{(\ell)}$ to be completely independent, for each $k$ and $\ell$, and conditioning on the new data vector $Z$; the random variable $T$, at (2.4), is seen to have variance equal to

$$\sum_{k=1}^{p} \mathrm{var}(I_{i_X}^{(k)} - J_{i_Y}^{(k)})(1 - 2K^{(k)})^2 = \sum_{k=1}^{p} \{\mathrm{var}(I_{i_X}^{(k)}) + \mathrm{var}(J_{i_Y}^{(k)})\}, \tag{2.10}$$

where the identity follows from the independence assumed earlier in this paragraph and from the fact that $1 - 2K^{(k)} = \pm 1$. Under the assumptions, $\mathrm{var}(I_{i_X}^{(k)}) = (EI_{i_X}^{(k)})(1 - EI_{i_X}^{(k)}) \leq E(I_{i_X}^{(k)})$, with an analogous result holding for $\mathrm{var}(I_{i_Y}^{(k)})$. Therefore, $S^2$, at (2.6), tends to overestimate the right-hand side of (2.10):

$$E(S^2) \geq \sum_{k=1}^{p} \{\mathrm{var}(I_{i_X}^{(k)}) + \mathrm{var}(J_{i_Y}^{(k)})\}.$$

This slight conservatism, and the log factor in the threshold $z_p$, provide opportunities for repairing errors that arise from failure of the independence assumption.

In the argument above, we defined $T$ to be the difference between the two sides of (2.2), rather than between the two sides of (2.3). Indeed, it can be awkward to estimate the variance of $T$ if we use (2.3) and do not have a good model for the distributions of $X^{(k)}$ and $Y^{(k)}$. There are ways of overcoming this difficulty, but we do not find them as attractive as working with (2.2).

An alternative approach to selecting $t$ could be based on standard cross-validation, taking $\theta$ to be the infimum of values $t$ that minimize the error-rate estimator,

$$\mathrm{CV}(t) = \frac{1}{m} \sum_{i=1}^{m} I\left(\min_{i_1 \neq i} \|X'_{i_1} - X'_i\| > \min_{1 \leq j \leq n} \|Y'_j - X'_i\|\right)$$
$$+ \frac{1}{n} \sum_{j=1}^{n} I\left(\min_{j_1 \neq j} \|Y'_{j_1} - Y'_j\| > \min_{1 \leq i \leq m} \|X'_i - Y'_j\|\right).$$

However, this technique has the disadvantage that it works only when $m$ and $n$ both exceed 1. Moreover, in most problems there is a continuum of values of $t$ that minimize $\mathrm{CV}(t)$, and so cross-validation does not give an explicit answer to the tuning-parameter choice problem.

2.3. *Example of mixed light- and heavy-tailed components.* When both light- and heavy-tailed data components are present in each data vector, and only a very small proportion of the components differ through perturbations, it can be particularly difficult to achieve good classification using standard distance-based methods, such as support vector machines. In the case of these approaches, the accumulation of noise from irrelevant components can drown out the signal in those few components that convey information for classification. Methods such as FDR, based on extrema, can bring substantial improvements in performance. However, when data distributions are heterogeneous, those techniques too can have difficulty.

To illustrate this point, assume for the sake of simplicity that all vector components are mutually independent. Suppose that $X$ consists of just $p^{1-\beta}$ components with standard normal distributions, where $\beta \in (0,1)$, and $p - p^{1-\beta}$ components having exponential distributions, for which $P(X^{(k)} > x) = e^{-x}$ when $x \geq 0$. Construct the $Y$ variable by adding $\mu = r \log p$, where $r > 0$, to just $p^{1-\beta}$ of the components of $X$, leaving the others unaltered.

If these special $p^{1-\beta}$ components are among those that have an exponential distribution, then we can write

(2.11)  $\quad \max_{1 \leq k \leq p} X^{(k)} = Q_1 + \log p + o_p(1),$

(2.12)  $\quad \max_{1 \leq k \leq p} Y^{(k)} = \max\{Q_1 + \log p, Q_2 + (1 - \beta + r) \log p\} + o_p(1),$

where $Q_1$ and $Q_2$ are asymptotically independent and have the extreme-value distribution function $\exp(-e^{-x})$. In both these expansions, we can consider $Q_1 + \log p$ to equal the maximum of the $p - 2p^{1-\beta}$ components of $X$ that have an exponential distribution and which exclude the $p^{1-\beta}$ components to which the perturbation $\mu$ is added to form $Y$, and $Q_2 + (1 - \beta + r)\log p$ to be the maximum of the $p^{1-\beta}$ components of $Y$ that are obtained by perturbing components of $X$.

It follows from (2.11) and (2.12) that if $r > \beta$, then $\max_k Y_k - \max_k X_k \to \infty$. More specifically, when $r > \beta$, the maximum of the components of a new vector $Z$ can be used to obtain asymptotically correct classification. This result does not hold if $r \leq \beta$.

On the other hand, if the perturbation $\mu$ is added to each of the $p^{1-\beta}$ components of $X$ that have a normal $N(0,1)$ distribution, and if none of the exponentially distributed components of $X$ is perturbed, then

$$\max_{1 \leq k \leq p} Y^{(k)} = \max[Q_1 + \log p, \{2(1 - \beta)(\log p)\}^{1/2} + r \log p] + o_p(1)$$

and $P(\max_k X_k = \max_k Y_k) \to 1$, unless $r \geq 1$. In particular, in this case, only when $r \geq 1$ is it possible to discriminate between the $X$ and $Y$ populations using extrema or FDR.

By way of contrast, we shall show in Section 3 that, no matter where the perturbations are added, the nearest-neighbor method produces asymptotically correct classification whenever $r > 2\beta - 1$. Since $1 > \beta > 2\beta - 1$, then the nearest-neighbor classifier enjoys greater sensitivity than the method based on extrema or FDR, no matter whether the perturbations are added to light- or heavy-tailed data components.

2.4. *Discussion of nearest-neighbor methods.* The versatility, performance and simplicity of NN classifiers are important factors in their popularity. As we show in this paper, NN methods also have significant potential for "robustification" and for fine-tuning through thresholding; both of these modifications lead to further improvements in performance. Nevertheless, well-known caveats about NN techniques should be mentioned.

Nearest-neighbor algorithms are most clearly suited to problems where the major departures among distributions are the results of differences in means, rather than differences in variances. To appreciate why NN classifiers can face challenges when the differences are principally in terms of variance, consider the elementary case where the variables $X^{(k)}$, for $1 \leq k \leq p$, are independent and identically distributed with zero mean and variance $\sigma_X^2$; the $Y^{(k)}$'s are likewise i.i.d., with zero mean and variance $\sigma_Y^2$; and $\sigma_X^2 < \sigma_Y^2$. If $Z$ comes from population $X$ then, as $p$ increases, the probability that the inequality

$$(2.13) \qquad \frac{1}{p}\sum_{k=1}^{p}(X^{(k)} - Z^{(k)})^2 < \frac{1}{p}\sum_{k=1}^{p}(Y^{(k)} - Z^{(k)})^2$$

holds tends to 1, since the left-hand side and right-hand side are, respectively, equal to $2\sigma_X^2 + o_p(1)$ and $\sigma_X^2 + \sigma_Y^2 + o_p(1)$. The probability that (2.13) holds when $Z$ is from population $Y$ also converges to 1, since in this setting the two sides of (2.13) are, respectively, $\sigma_X^2 + \sigma_Y^2 + o_p(1)$ and $2\sigma_Y^2 + o_p(1)$. Therefore, no matter what population $Z$ is from, the simple NN classifier will, with probability converging to 1 as $p \to \infty$, assign $Z$ to the population with smaller variance, that is, to population $X$. This will hold true for samples of any sizes $m$ and $n$, provided those quantities are kept fixed as $p$ diverges.

The result is quite different if the two populations have equal variances but unequal means. There, the probability that $Z$ is correctly allocated by a NN classifier typically converges to 1 as $p \to \infty$, if there are sufficiently many sufficiently large differences among means. Although in Section 3 we shall permit distributions to take very different forms among components, the differences with real leverage for classification will be those among means. The process of thresholding, which converts continuous measurements into zero–one data, tends to remove problems caused by differences among variances, although to some extent it converts differences among means into differences among variances; recall that a zero–one variable with mean $q$ has variance $q(1-q)$. However, as we shall show in Section 3, this does not cause significant difficulty.

## 3. Theoretical properties.

3.1. *Summary.* The models that we shall use to describe the $X$ and $Y$ vectors will differ through perturbations (location changes), $\mu^{(k)}$, added to individual components. The models will be constructed so as to admit considerable heterogeneity among the distributions, as well as to allow dependence; see Section 3.6 for discussion of the latter. In Sections 3.2 and 2.3 we shall describe the density, size and scalability of the perturbations and marginal distributions. Classification boundaries will be discussed in Section 3.4. The principles introduced there will dictate the context of the main theoretical results given in Sections 3.5 and 3.6. These results will reflect difficult classification problems, where configurations are close to optimal classification boundaries. Our main theorems will be stated under the assumption that the number of dimensions, $p$, diverges, while the sample sizes, $m$ and $n$, are held fixed.

3.2. *Relationship between marginal distributions of $X$ and $Y$.* For sequences $b_p$ and $c_p$ depending on $p$, we write $b_p \asymp c_p$ to mean that the ratio $b_p/c_p$ is bounded above zero and below infinity, as $p$ diverges. Given a sequence $a_p$ diverging to infinity, and a constant $\beta \in (\frac{1}{2}, 1)$, we shall say that

(3.1) the sequence $\mu^{(1)}, \ldots, \mu^{(p)}$ "has asymptotic density $p^{-\beta}$ and is on the scale $a_p$," if (i) the number, $N_p$ say, of nonzero $\mu^{(k)}$'s satisfies $N_p \asymp p^{1-\beta}$ and (ii) none of the nonzero $\mu^{(k)}$'s is less than $a_p$.

The perturbations $\mu^{(k)}$ will be added to the respective components of $X$ to create a vector with the distribution of $Y$. Therefore, our model for the way in which the marginal distributions of $X$ and $Y$ are related will be that

(3.2) for $1 \leq k \leq p$, $Y^{(k)}$ is distributed as $X^{(k)} + \mu^{(k)}$, where the sequence $\mu^{(1)}, \ldots, \mu^{(p)}$ has asymptotic density $p^{-\beta}$ and is on the scale $a_p$, with $\beta \in (\frac{1}{2}, 1)$.

Condition (3.2) relates only to the number of $\mu^{(k)}$'s that are different from zero, not to the order of the nonzero values in the sequence $\mu^{(1)}, \ldots, \mu^{(p)}$. In particular, the assumption is much less stringent than it would be if it were supposed that the indices of the nonzero $\mu^{(k)}$'s were distributed according to a particular random process. The latter constraint is implicit whenever mixture models are assumed.

In (3.1) and (3.2), we choose $\beta \in (\frac{1}{2}, 1)$ since classification in the case $\beta \leq \frac{1}{2}$ is relatively easy (indeed, root-$n$ consistent estimation is generally possible when $\beta < \frac{1}{2}$), and since nontrivial, asymptotically correct classification is impossible if $\beta = 1$.

3.3. *Scalability.* We shall use the phrase, "the marginal distributions of $X$ are continuous and scalable," to mean that, for each $r \in (0, 1)$, the equation

$$\text{(3.3)} \qquad \sum_{k=1}^{p} P(X^{(k)} > a_p) = p^{1-r}$$

has a unique solution $a_p = a_p(r)$, and that for each $\varepsilon \in (0, r)$ there exists $C = C(\varepsilon) \in (0, 1)$ such that, for all sufficiently large $p$,

$$\text{(3.4)} \qquad \sum_{k=1}^{p} P(X^{(k)} > Ca_p) \leq p^{1-r+\varepsilon}.$$

In particular, if the $X^{(k)}$'s are identically distributed as $X^{(0)}$, then the common distribution is scalable if, when $a_p$ is defined by $P(X^{(0)} > a_p) = p^{-r}$, for each $\varepsilon \in (0, r)$ there exists $C \in (0, 1)$, such that $P(X^{(0)} > Ca_p) \leq p^{1-r+\varepsilon}$. Scalable distributions include the normal, and other exponentially decreasing distributions such as the Subbotin, with probability density function $f$ given by

$$\text{(3.5)} \qquad f(x) = C_\gamma^{-1} \exp(-|x|^\gamma/\gamma),$$

where $\gamma > 0$ and $C_\gamma = 2\Gamma(1/\gamma)\gamma^{(1/\gamma)-1}$. See Donoho and Jin [5] for an account of the interest in, and applications of, the Subbotin distribution. Scalable distributions also include regularly varying distributions such as the Pareto, for which

$$\text{(3.6)} \qquad P(X^{(k)} > x) = x^{-\gamma},$$

when $x > 1$, where $\gamma > 0$. Nonscalable distributions have extremely light upper tails, for example, the extreme-value distribution for which $P(X > x) = \exp(-e^x)$.

Of course, scalability of the marginal distributions of $X$ does not require the $X^{(k)}$'s to be identically distributed. A particularly simple, nonidentically distributed example is that where $N_1(p)$ of the components $X^{(k)}$ are distributed as $X^{(0)}$, say; the other $N_2(p) = p - N_1(p)$ components have distribution functions that dominate that of $X^{(0)}$, in the sense that $P(X^{(k)} \leq t) \geq P(X^{(0)} \leq t)$ for all $1 \leq k \leq p$ and all $t \geq t_0$, say; the distribution of $X^{(0)}$ is scalable, in the sense described in the previous paragraph; and $N_1(p) \sim p$ as $n \to \infty$. This model, and Theorems 1 and 2 below, permit a rigorous account of performance of the nearest-neighbor classifier in the context of the examples discussed in Section 2.3.

3.4. *Detection and classification boundaries.* In this subsection, we assume that all the marginal distributions of $X$ are identical to that of $X^{(0)}$, say, and we take each of the $p^{1-\beta}$ nonzero values of $\mu^{(k)}$ to equal $a_p$, defined by $P(X^{(0)} > a_p) = p^{-r}$, where $\beta \in (\frac{1}{2}, 1)$ and $r \in (0, 1)$. Theorems 1 and 2, below, imply that in this case the robust nearest-neighbor classifier defined by (2.9) will asymptotically correctly classify data, provided that

$$(3.7) \qquad 1 - 2\beta + r > 0.$$

That is, if (3.7) holds, and even when $m = n = 1$ (i.e., when there is only one training data value from each population), the probability that the classifier at (2.9) correctly assigns $Z$, no matter whether it comes from the $X$ or the $Y$ population, converges to 1 as $p \to \infty$.

Conversely, if $(\beta, r)$ lies strictly below the boundary described by the line

$$(3.8) \qquad 1 - 2\beta + r = 0,$$

then the probability of correct classification fails to converge to 1. Moreover, the same boundary plays the same role (i.e., as the border that separates classifiable and nonclassifiable cases) if we use a truncated standard nearest-neighbor method. The latter technique requires the data distributions to have several finite moments, whereas the approach suggested in our paper is far more robust than conventionally truncated nearest-neighbor methods. It is significant that the boundaries are identical in the cases of robust and nonrobust nearest-neighbor methods. In particular, the greater robustness of our approach does not come at the price of reduced effectiveness.

To define standard truncated nearest-neighbor classifiers, let $X_{ij}^{\text{tr}} = X_{ij} I(X_{ij} > t)$, $Y_{ij}^{\text{tr}} = Y_{ij} I(Y_{ij} > t)$ and $Z_j^{\text{tr}} = Z_j I(Z_j > t)$, respectively, where $t$ denotes the truncation point. The corresponding truncated vectors are $X_i^{\text{tr}} = (X_{ij}^{\text{tr}})$, $Y_i^{\text{tr}} = (Y_{ij}^{\text{tr}})$ and $Z^{\text{tr}} = (Z_j^{\text{tr}})$. We apply the standard nearest-neighbor classifier to the truncated datasets $\{X_1^{\text{tr}}, \ldots, X_m^{\text{tr}}\}$ and $\{Y_1^{\text{tr}}, \ldots, Y_n^{\text{tr}}\}$, instead of to the original data. That is, we assign $Z$ to the $X$ population if $Z^{\text{tr}}$ is nearer to at least one of $X_i^{\text{tr}}$'s than it is to any of the $Y_i^{\text{tr}}$'s, and we assign it to the $Y$ population otherwise. Assume that the random variables $X_{i_1 j_1} - E(X_{i_1 j_1})$ and $Y_{i_2 j_2} - E(Y_{i_2 j_2})$ are all independent and identically distributed, with the distribution of $U$, say, and that the scalability condition holds. It can be proved that if $q = p^{1-\beta}$; if the truncation point $t$ does not exceed $\nu$; if $\nu = a_p$, where $a_p$ satisfies (3.3) [or equivalently, $P(U > a_p) = p^{-r}$]; and if $(\beta, r)$ lies strictly below the boundary given by (3.8); then $pE\{U^4 I(U > t)\}/(q\nu^2)^2$ is bounded away from zero. Moreover, it is shown by Hall, Pittelkow and Ghosh [8] that if, along a subsequence of values of $p$, $pE\{U^4 I(U > t)\}/(q\nu^2)^2$ does not converge to zero, then the probability of correct classification fails to converge to 1. Similarly, if $(\beta, r)$ lies above the boundary, then the probability of correct classification converges to 1. This establishes the implications of the boundary in the case of standard truncated nearest-neighbor classifiers, and its implications for our truncated form are similar.

In some problems, and when $m = n = 1$, the boundary at (3.8) is identical to that for an optimal classifier, implying that the robust nearest-neighbor approach has asymptotically optimal performance. However, the classifiers for which this boundary is known require the marginal distribution to be known; our truncated, thresholded nearest-neighbor approach is not subject to that requirement.

For example, in the Subbotin case represented by (3.5), with $0 < \gamma \leq 1$; and in the Pareto case given by (3.6), when $\gamma > 0$; it is known [5, 14] that the boundary represented by (3.8) is the optimal boundary for signal detection. It can be proved from this result that it is also the optimal boundary for classification, when $m = n = 1$. In the Subbotin case where $\gamma > 1$, alternative methods, such as Donoho and Jin's [5] higher-criticism method and the approaches suggested by Ingster [11, 12, 13], give a lower optimal boundary even when $m = n = 1$ and, hence, permit classification in cases where robust nearest-neighbor methods do not.

3.5. *Case of independent components.* The error rates of the classifier at (2.8) are defined to be the probability that $Z$ is misclassified as coming from $Y$ when it is really from $X$ and the probability of misclassification of $Z$ as coming from $X$ when it is actually from $Y$.

Note that, if $t$ is sufficiently small, then it is possible to have $P(X^{(k)} > t) = P(Y^{(k)} > t) = 1$, uniformly in $1 \leq k \leq p$ and in $p$. In this case, the ratio $T(t)/S(t)$ is not well defined. To remove pathologies such as this, we modify the definition of $\theta$, at (2.7), by insisting that, for some fixed $t_0$ sufficiently large, only values $t \geq t_0$ be considered. In the theorem below, we hold $m$ and $n$ fixed and let $p$ increase without bound.

THEOREM 1. *If the components of $X$ are independent, and the components of $Y$ are independent; if the marginal distributions are related by (3.2), and are continuous and scalable; if, for $r \in (0,1)$, the quantity $a_p = a_p(r)$, defined by (3.3), diverges to infinity but at a rate no faster than $p^D$ for some $D > 0$, as $p$ increases; if the pair $(\beta, r)$ is above the classification boundary, in the sense that (3.7) holds; and if $z_p$ is given by (2.8), where $\xi_p$ diverges more slowly than $p^\varepsilon$ for each $\varepsilon > 0$; then, as $p \to \infty$ for fixed $m$ and $n$, the error rates of the classifier at (2.9) converge to zero.*

The assumption in Theorem 1 that $a_p = O(p^D)$, for some $D > 0$, is satisfied if, for example, $\sup_k E(|X^{(k)}|^\varepsilon) < \infty$ for some $\varepsilon > 0$.

3.6. *Case of dependent components.* As in (3.2), we take the distributions of the components of $Y$ to be translations of those of the respective components of $X$. In particular, given stochastic processes $U_1, \ldots, U_p$ and $U_1^\#, \ldots, U_p^\#$, each with the same $p$-variate distribution, we define

$$(3.9) \qquad X^{(k)} = U_k + \nu_k, \qquad Y^{(k)} = U_k^\# + \nu_k + \mu^{(k)}.$$

The challenge is to model the degree of dependence among marginals and, at the same time, to permit the marginal distributions to vary in shape, as well as location, from one component to another. This is done through an exponentiated moving average process, defined in part (a) of (3.10):

(3.10) (a) $U_k = \sum_{j \geq 1} \omega_j W_{j+k}^{\alpha_k}$, where the nonnegative random variables $W_j$ are independent and identically distributed as $W$; (b) for all $w$, $P(W \leq w) < 1$; (c) for some $c > 0$, $E(W^c) < \infty$; (d) the distribution of $W$ has a bounded probability density; (e) the constants $\alpha_k$ are permitted to be functions of $p$ as well as $k$, and for some $C > 1$, for all $p$ and for all $1 \leq k \leq p$, $C^{-1} \leq \alpha_k \leq C$; (f) for some $C > 0$, for some $\omega \in (0,1)$ and for all $j \geq 1$, $|\omega_j| \leq C\omega^j$; and (g) at least one $\omega_j$ is strictly positive.

The $\nu_k$'s are taken to be uniformly bounded, and the $\mu^{(k)}$'s to have properties similar to those at (3.2):

(3.11) (a) $\nu_k$ and $\mu^{(k)}$ are functions of $p$ as well as $k$; (b) for a fixed constant $C > 0$, $|\nu_k| \leq C$ for all $p$ and for all $1 \leq k \leq p$; and (c) given $r \in (0,1)$ and $\beta \in (\frac{1}{2}, 1)$ and with $a_p$ defined by (3.3), the sequence $\mu^{(1)}, \ldots, \mu^{(k)}$ has asymptotic density $p^{-\beta}$ and is on the scale $a_p$.

The "continuity" part of the assumption, in Theorem 1, that the marginal distributions of $X$ are continuous and scalable, is taken care of by (3.10)(d). However, we also need scalability, as well as a version of that condition in the case of logarithmically spaced marginals. For the latter, (3.12) is sufficient: defining $\pi_k(t) = P(X^{(k)} \geq t)$, we ask that

(3.12) for each $B, \varepsilon > 0$, there exists $t' = t'(B, \varepsilon)$ such that, if $\ell_p$ denotes the integer part of $B \log p$, then $p^\varepsilon \sum_{0 \leq k \leq (p-h)/\ell_p} \pi_{k\ell_p + h}(t) \geq \sum_{1 \leq k \leq p} \pi_k(t)$ for all $0 \leq h \leq \ell_p$ and all $t \geq t'$.

In assumption (3.10), parts (a) and (e) imply that the $U_k$ process is a generalized moving average with geometrically decaying coefficients. The generalization, through raising $W_{j+k}$ to the power $\alpha_k$, allows the distribution of $U_k$ to be varied substantially from one component to another. In particular, the tail weights can be very different; smaller $\alpha_k$'s give distributions with lighter tails.

To interpret parts (b) and (g) of (3.10), note that if $P(U_k \leq C) = 1$ for some $C > 0$ and for all $k$, then the problem of discriminating between $X$ and $Y$, on the basis of location shifts to the right, is relatively simple. Part (b), which asserts that the upper tail of the distribution of $W$ is unbounded, together with (g), which asks that at least one contribution $\omega_j W_{j+k}^{\alpha_k}$ to $U_k$ be positive, permit us to avoid this degeneracy. Part (c) of (3.10) is a very weak moment assumption and, in particular, permits the distribution of $W$ to be so heavy tailed that it lies in the domain of attraction of a stable law.

In (3.11), parts (a) and (b) permit the $\nu_k$'s to vary quite generally, subject only to being bounded. Condition (3.12) holds true trivially if the marginal distributions are all identical and can be shown to be valid under other heterogeneous models.

Theorem 2, below, is a version of Theorem 1 for dependent data. As in the case of Theorem 1, we modify the definition of $\theta$, at (2.7), by considering only values $t \geq t_0$, for $t_0$ fixed but sufficiently large.

THEOREM 2. *If the joint distributions of the components of $X$ and $Y$ are given by (3.9), with the quantities there generated as described by (3.10) and (3.11); if the marginal distributions of $X^{(k)}$ are scalable, and satisfy (3.12); if, for $r \in (0, 1)$, the quantity $a_p = a_p(r)$, defined by (3.3), diverges no faster than $p^D$ for some $D > 0$, as $p$ increases; if the pair $(\beta, r)$ lies above the classification boundary, in the sense that (3.7) holds; and if $z_p$ is given by (2.5), where $\xi_p$ diverges more slowly than $p^\varepsilon$ for each $\varepsilon > 0$; then, as $p \to \infty$ for fixed $m$ and $n$, the error rates of the classifier at (2.9) converge to zero.*

## 4. Numerical properties.

4.1. *Microarray data.* As a practical example, we compared the performance of the thresholded method with the nearest-neighbor method on the BRCA dataset [6, 9], which we obtained from http://www.nejm.org/general/

content/supplemental/hedenfalk/index.html. This dataset contains microarray data from patients with breast cancer, caused by two different types of mutations, labelled BRCA1 and BRCA2. The expression level of each of 3226 genes was measured in each patient, and there are 7 patients with BRCA1 and 8 patients with BRCA2.

This dataset (and indeed many gene microarray datasets) is very suited to our thresholded method. For a start, it is a dataset with very high dimension and low sample size. Furthermore, it is expected that only a few genes will be differentially expressed between the two types of cancer, so the difference between the populations is sparse. Lastly, the underlying distributions of the gene expressions are likely to be both heavy-tailed and with significant dependence among genes, which nearest-neighbor traditionally does poorly at, especially in comparison with the thresholded method.

We tested the two methods on this dataset by calculating the cross-validation performance, where we classify each patient according to all the other patients and calculate the classification rate. For the nearest-neighbor method, cross-validation correctly classified 11 out of the 15 patients. Our thresholded method did a lot better; with $z_p = 0.5(\ln p)^{1/2}$, all 15 patients were classified correctly under cross-validation. In fact, this happened when we set the coefficient of $(\ln p)^{1/2}$ in $z_p$ to be anywhere between 0.35 and 0.5.

4.2. *Simulated data.* As an additional test, we also compared the thresholded method with the nearest-neighbor method for simulated data. We compare the two methods in the area of the $\beta$–$r$ plane where classification is possible ($r > 2\beta - 1$), but not easy ($\beta < \frac{1}{2}$ or $r > 1$). Overall, we found that in cases where standard nearest-neighbor does not perform well, the thresholded method improves on it. We look at some of these cases.

4.2.1. *Independent heavy-tailed marginal distributions.* Nearest-neighbor methods do not do very well when the marginal distributions of the components of $X$ (and $Y$) are heavy-tailed (i.e., go to 0 slower than a normal distribution). We compared the methods for simple models where $m = n = 1$ and each of the components of $X$ are independent and have identical Student's-$t$ distributions. By varying the degrees of freedom, we can observe the behavior of the methods relative to the heaviness of the tails.

For this case, if we are given a threshold $t$, the success rate of the algorithm can be approximated very accurately, for any $\beta$ and $r$, by looking at the contribution of each dimension to $T(t)$. By varying $t$ we can calculate the optimal threshold, which we call the a priori optimal threshold, and also the best possible performance of the classifier. However, we are not usually given the threshold, so this is an upper limit on the success rate. Instead, we compare the classifiers with empirically chosen thresholds, on simulated data with $p$ up to 20,000.

We found that, for sufficiently heavy tails, the thresholded method dominates standard nearest-neighbor in all areas of the $\beta$–$r$ plane. In fact, the success rate of the thresholded method actually improves for heavier tails. As the tails get lighter (the d.f. gets larger), the success rate declines, and nearest-neighbor does better in a small area in the plane, which grows and moves around as the tail weight decreases. For small d.f., this area occurs at high $\beta$ and $r$ [see Figure 1(a)]; for larger d.f., this area occurs at low $\beta$ and $r$ neither high nor low [see Figure 1(b)].

The thresholded method also dominates nearest-neighbor if we use the a priori optimal thresholds, for sufficiently heavy tails. If the tails are not heavy enough, nearest-neighbor works better for low $\beta$ and $r$.

We found that the best performance of the thresholded method is achieved when we take $z_p$ in (2.7) to be $c(\ln p)^{1/2}$, where $c$ is a constant. The value of $c$, which maximizes the success rate, lies between 0.3 and 0.9, depending on $\beta$ and $r$. However, the best success rate achieved with an empirically chosen threshold is worse than that achieved with the a priori threshold, because the empirical threshold is not constant for constant $z_p$. Figure 2 estimates the distribution of the chosen threshold for various cases when $z_p$ is close to optimal. Figure 3 shows how the value of the threshold affects the success rate, while Figure 4 shows how the value chosen for $z_p$ affects the success rate. In both of these figures, the curves represent the thresholded method, while the horizontal lines show the performance of the nearest-neighbor method for comparison.

4.2.2. *Dependent normal marginal distributions.* Another case where standard nearest-neighbor methods perform badly is when the components of $X$ are dependent on each other. We compared the methods for varying degrees and types of dependence; for example, when the components of $X$ are moving averages of independent standard normal variables, or weighted moving averages, or an autoregressive process $X^{(i+1)} = \alpha X^{(i)} + (1-\alpha)N^{(i)}$, where $N$ is a sequence of independent standard normal variables.

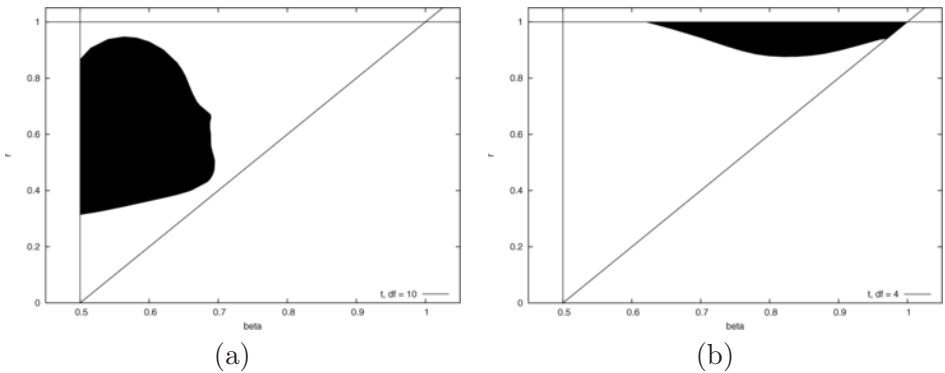

Fig. 1. *Areas where the two methods perform better, for heavy-tailed distributions. The nearest-neighbor method performs better in the shaded area; otherwise, the thresholded method is better.* (a) *t distributions, d.f.* = 4, (b) *t distributions, d.f.* = 10.

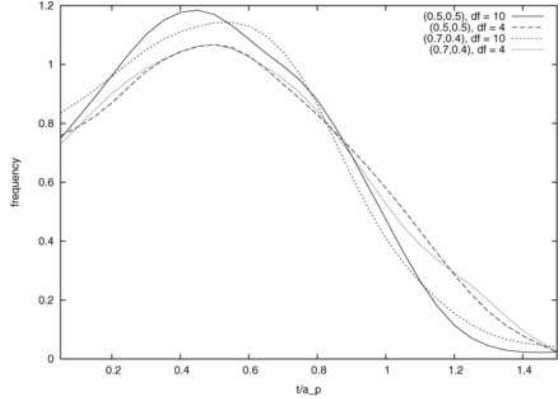

FIG. 2. *Estimated distribution of thresholds produced with t distributions at $z_p = 0.55(\ln p)^{1/2}, p = 20{,}000$, at various $(\beta, r)$ and degrees of freedom.*

Again, we found that for sufficient levels of dependence, the thresholded method dominates nearest-neighbor for all $(\beta, r)$. For weaker levels of dependence, the nearest-neighbor method works better in a small area at small $\beta$ and $r$ neither small nor large (see Figure 5), and this region grows with decreasing dependence. We found that the strength of the dependence [e.g., $\text{cov}(X^{(i)}, X^{(i+1)})$] affects the size of this region more than the length of the dependence (the number of components of $X$ dependent on a given component).

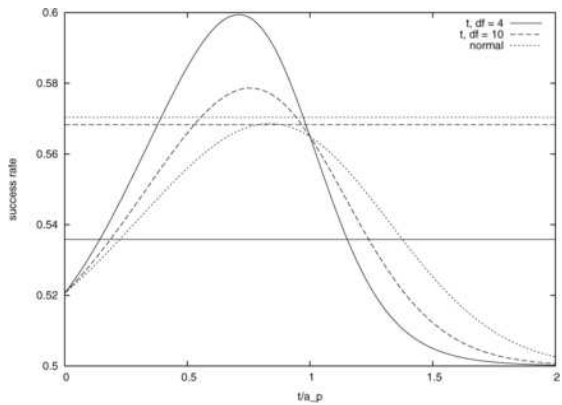

FIG. 3. *Success rate vs. threshold (as a proportion of shift amount) for $p = 20{,}000, (\beta, r) = (0.7, 0.4)$.*

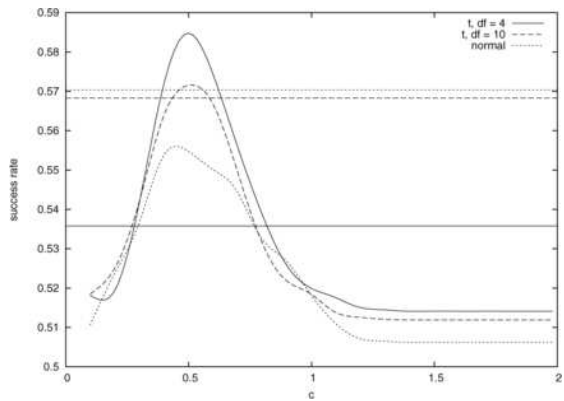

FIG. 4. *Success rate vs. c for $p = 20{,}000$, $(\beta, r) = (0.7, 0.4)$, where $z_p = c(\ln p)^{1/2}$.*

As with the heavy-tailed case, taking $z_p = c(\log p)^{1/2}$ optimizes the success rate, with $c$ taking similar values as before. However, the overall success rate of the thresholded method is worse than for an equivalent independent case. The behavior of the chosen threshold, and its effect on the success rate, is similar to its behavior for heavy-tailed distributions.

4.2.3. *Independent normal marginal distributions.* For comparison, we also looked at the case where the components of $X$ were independent and normally distributed. Here, the thresholded method does not dominate nearest-neighbor, which works better for low $\beta$ (approximately $\beta < 0.65$). This is consistent with heavy-tailed distributions as the tails get lighter. The behavior of the chosen threshold, and its effect on the success rate, is again similar to its behavior for heavy-tailed distributions. The overall success rate is worse than for heavy-tailed distributions, but better than that for dependent distributions.

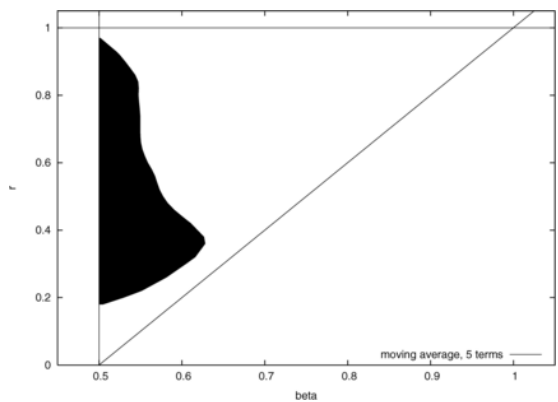

FIG. 5. *Areas where the two methods work best, for moving averages of 5 normal random variables.*

4.2.4. *Larger samples.* The above scenarios all involved $m = n = 1$. As the sample sizes $m$ and $n$ increase, but are kept equal, the classification success rate of both methods increase. As $m$ and $n$ increase, the thresholded method outperforms the nearest-neighbor method for a greater range of the $\beta$–$r$ plane, although the difference is slight up to $m = n = 10$ (the upper limit of our testing).

When the sample sizes are not equal, the thresholded method performs better when $m$ is smaller, if $m + n$ is kept constant. In fact, although increasing $m$ or $n$ while keeping the other fixed generally increases the classification rate, it is possible to decrease the classification rate by increasing $m$ while keeping $n$ fixed (e.g., when $n = 1$). As the effectiveness of the nearest-neighbor method stays largely the same, the thresholded method outperforms the nearest-neighbor method for much larger areas of the $\beta$–$r$ plane, when $m < n$, and is much less effective for $m > n$.

## REFERENCES

[1] COVER, T. M. (1968). Rates of convergence for nearest neighbor procedures. In *Proceedings of the Hawaii International Conference on System Sciences* (B. K. Kinariwala and F. F. Kuo, eds.) 413–415. Univ. Hawaii Press, Honolulu.
[2] COVER, T. M. and HART, P. E. (1967). Nearest neighbor pattern classification. *IEEE Trans. Inform. Theory* **13** 21–27.
[3] DASARATHY, B. V. (1991). *Nearest Neighbor (NN) Norms: NN Pattern Classification Techniques.* IEEE Computer Society, Los Alamitos, CA.
[4] DEVROYE, L., GYÖRFI, L. and LUGOSI, G. (1996). *A Probabilistic Theory of Pattern Recognition.* Springer, New York. MR1383093
[5] DONOHO, D. J. and JIN, J. (2004). Higher criticism for detecting sparse heterogeneous mixtures. *Ann. Statist.* **32** 962–994. MR2065195
[6] EFRON, B. (2004). Large-scale simultaneous hypothesis testing: The choice of a null hypothesis. *J. Amer. Statist. Assoc.* **99** 96–104. MR2054289
[7] FRITZ, J. (1975). Distribution-free exponential error bound for nearest neighbor pattern classification. *IEEE Trans. Inform. Theory* **21** 552–557. MR0395379
[8] HALL, P., PITTELKOW, Y. and GHOSH, M. (2008). Theoretical measures of relative performance of classifiers for high-dimensional data with small sample sizes. *J. Roy. Statist. Soc. Ser. B* **70** 159–173.
[9] HEDENFALK, I., DUGGAN, D., CHEN, Y., RADMACHER, M., BITTNER, M., SIMON, R., MELTZER, P., GUSTERSON, B., ESTELLER, M., RAFFELD, M., YAKHINI, Z., BEN-DOR, A., DOUGHERTY, E., KONONEN, J., BUBENDORF, L., FEHRLE, W., PITTALUGA, S., GRUVBERGER, S., LOMAN, N., JOHANNSSON, O., OLSSON, H., WILFOND, B., SAUTER, G., KALLIONIEMI, O.-P., BORG, A. and TRENT, J. (2001). Gene expression profiles in hereditary breast cancer. *N. Engl. J. Med.* **344** 539–548.
[10] HOLST, M. and IRLE, A. (2001). Nearest neighbor classification with dependent training sequences. *Ann. Statist.* **29** 1424–1442. MR1873337

[11] INGSTER, Y. I. (1999). Minimax detection of a signal for $l^n$-balls. *Math. Methods Statist.* **7** 401–428. MR1680087
[12] INGSTER, Y. I. (2001). Adaptive detection of a signal of growing dimension. I. Meeting on mathematical statistics. *Math. Methods Statist.* **10** 395–421. MR1887340
[13] INGSTER, Y. I. (2002). Adaptive detection of a signal of growing dimension. II. *Math. Methods Statist.* **11** 37–68. MR1900973
[14] JIN, J. (2002). Detection boundary for sparse mixtures. Unpublished manuscript.
[15] KULKARNI, S. R. and POSNER, S. E. (1995). Rates of convergence of nearest neighbor estimation under arbitrary sampling. *IEEE Trans. Inform. Theory* **41** 1028–1039. MR1366756
[16] PSALTIS, D., SNAPP, R. R. and VENKATESH, S. S. (1994). On the finite sample performance of the nearest neighbor classifier. *IEEE Trans. Inform. Theory* **40** 820–837.
[17] SHAKHNAROVICH, G., DARRELL, T. and INDYK, P. (2006). *Nearest-Neighbor Methods in Learning and Vision*. MIT Press, Boston.
[18] WAGNER, T. J. (1971). Convergence of the nearest neighbor rule. *IEEE Trans. Inform. Theory* **17** 566–571. MR0298829

DEPARTMENT OF MATHEMATICS AND STATISTICS
UNIVERSITY OF MELBOURNE
PARKVILLE, VIC 3010
AUSTRALIA
E-MAIL: y.chan@ms.unimelb.edu.au
        P.Hall@ms.unimelb.edu.au